\documentclass[a4paper,fleqn]{cas-sc}

\usepackage[authoryear,longnamesfirst]{natbib}

\usepackage{graphicx}%
\usepackage{multirow}%
\usepackage{amsmath,amssymb,amsfonts}%
\usepackage{amsthm}%
\usepackage{mathrsfs}%
\usepackage[title]{appendix}%
\usepackage{xcolor}%
\usepackage{textcomp}%
\usepackage{manyfoot}%
\usepackage{booktabs}%
\usepackage{algorithm}%
\usepackage{algorithmicx}%
\usepackage{algpseudocode}%
\usepackage{listings}%
\usepackage{placeins}
\usepackage{lineno}
\usepackage{setspace}

\newtheorem{theorem}{Theorem}
\newtheorem{remark}{Remark}%
\newtheorem{lemma}{Lemma}%

\def\tsc#1{\csdef{#1}{\textsc{\lowercase{#1}}\xspace}}
\tsc{WGM}
\tsc{QE}

\begin{document}

\let\WriteBookmarks\relax
\def\floatpagepagefraction{1}
\def\textpagefraction{.001}

\shorttitle{}    

\shortauthors{}  

\title [mode = title]{Solution to an infinite horizon McKean-Vlasov FBS$\Delta$E and related LQ control problems with law-delay}  

\tnotemark[1] 

\tnotetext[1]{} 

%

\author[concordia]{Duocheng Wang}
\ead{duocheng.wang@concordia.ca}

\address[concordia]{
Department of Mathematics and Statistics,
Concordia University,
Montreal, QC H3G 1M8, Canada
}

\address[concordia]{
Department of Mathematics and Statistics,
Concordia University,
Montreal, QC H3G 1M8, Canada
}


\begin{abstract}
Motivated by infinite horizon stochastic optimal control problems with law delay, this paper investigates a class of fully coupled McKean–Vlasov forward-backward stochastic difference equations (FBS$\Delta$Es) with delayed and anticipated law dependence. The main challenges arise from the infinite time horizon, the fully coupled structure, and the interaction between discounting and time-shifted distributional terms. By introducing a new monotonicity condition, we establish the existence and uniqueness of adapted solutions through the continuation method.  As an application, the infinite horizon linear-quadratic McKean–Vlasov control problem with law delay is studied and the unique optimal control is obtained through the associated FBS$\Delta$E.
\end{abstract}




\maketitle
\doublespacing


\section{Introduction}

Forward-backward stochastic differential equations (FBSDEs) have been widely studied in stochastic analysis and optimal control. The well-posedness of fully coupled FBSDEs was established through the four-step scheme \cite{ma1994solving} and the continuation method \cite{hu1995solution,peng1999fully}. Meanwhile, McKean--Vlasov stochastic systems have attracted considerable attention due to their applications in mean-field control and large population systems. Stochastic maximum principles for such systems have been developed in \cite{andersson2011maximum,buckdahn2016stochastic,nie2022extended,dong2022maximum,buckdahn2026global}. Related problems involving state or distribution delays have also been investigated; see \cite{chen2019maximum,zhang2021stochastic,guo2024stochastic}.

Discrete-time stochastic control problem is also important. Following \cite{lin2014maximum}, various stochastic maximum principles have been established, such as \cite{wu2022maximum,dong2023maximum,niu2026stochastic}. For infinite-horizon systems, \cite{ji2024infinite} studied backward stochastic difference equations and related recursive control problems.  \cite{song2023solvability} investigated fully coupled forward-backward stochastic difference equations (FBS$\Delta$Es). However, as far as we know, the solvability of infinite horizon fully coupled McKean-Vlasov FBS$\Delta$Es involving both delayed and anticipated law dependence has not been investigated.

Compared with finite horizon problems, the infinite horizon setting presents additional challenges arising from the absence of a prescribed terminal condition, the need for weighted square integrability, and the interaction between exponential weights and delayed or anticipated law dependence.

Motivated by these challenges, this paper focuses on infinite-horizon fully coupled McKean--Vlasov FBS$\Delta$Es with delayed and anticipated law dependence. Under suitable Lipschitz and monotonicity conditions, we establish the existence and uniqueness of adapted solutions through the continuation method. A new monotonicity condition is proposed to account for the effects of delayed control distributions, and weighted index-shifting techniques are employed to handle the resulting anticipated terms. As an application, an infinite horizon linear quadratic (LQ) McKean--Vlasov control problem with law delay is considered and its unique optimal control is obtained through the associated FBS$\Delta$E.

The remainder of this paper is organized as follows. In section 2, the well-posedness of the infinite horizon FBS$\Delta$E is established. In section 3,  the corresponding linear quadratic optimal control problem is investigated.

\section{The forward-backward stochastic difference equations
}
Let $(\Omega,\mathcal{F},\mathbb{P})$ be a probability space.  Let $W$ be a $\mathbb{R}$-valued square integrable martingale
process with independent increments and the predictable representation property, such that $\mathbb{E}\left[\Delta W_k\Delta W_l\right]={\delta_{kl}}$, where $\Delta W_k=W_{k+1}-W_k$ and $\delta$ is the Dirac function. Let $\mathcal{F}_0\subset \mathcal{F}$ be a sub $\sigma$-algebra independent of $W$, and  {$\mathbb{F}=(\mathcal{F}_k)_{0\le k\le N}$} be the filtration defined by $\mathcal{F}_k=\mathcal{F}_0\lor\sigma(\Delta W_0,...,\Delta W_{k-1})$. We define the following notations:

$$
\begin{aligned}
C[-\tau,0]^n:= & \left\{\varphi_k:[-\tau,0] \rightarrow \mathbb{R}^n \text { is }\mathcal{F}_0 \text {-measurable and } \mathbb{E}\sum_{k=-\tau}^0\|\varphi_k\|^2<+\infty\right\} \\
L_{\mathcal{F}}^{2,\lambda}\left(a, b ; \mathbb{R}^n\right):= & \left\{\varphi_k, a \leq k \leq b, \text { is an } \mathcal{G}_k \text {-adapted stochastic process s.t. } \mathbb{E} \sum_{k=a}^b(1\land e^{-\lambda k})\|\varphi_k\|^2<+\infty\right\}
\end{aligned}
$$
with $\mathcal{G}_k=\mathcal{F}_{k\lor 0}$.  {Here $\tau>0$ is an integer-valued constant time delay parameter} and $\lambda>0$ is the discount factor.

Denote $\mathcal{P}_2(\mathbb{R}^n)$ the space of square-integrable probability measures on $(\mathbb{R}^n,\mathcal{B}(\mathbb{R}^n))$, $\mathcal{T}=\{0,1,...\}$, $\Delta\psi_k=\psi_{k+1}-\psi_k$ for all adapted processes $\psi$. Let $\mathcal{L}(X)$ be the probability distribution of $X$. We introduce the following McKean-Vlasov  {FBS$\Delta$E} contains delay and anticipated term:
\begin{equation}\label{fbsde}
\begin{cases}\Delta x_k=b\left(k, u_{1,k}\right)+\sigma\left(k,u_{1,k}\right)\Delta W_k, & k\in\mathcal{T}, 
\\ 
y_k+z_k\Delta W_k=e^{-\lambda}\left[y_{k+1}+f(k+1,u_{1,k+1})\right], & k\in\mathcal{T}, 
\\ 
x_k=\varphi_k, &  {k \in[-\tau, 0]}, \\ 
\lim_{N\to+\infty}e^{-\lambda N}\mathbb{E}\left[\|y_N\|^2+\|z_N\|^2\right]=0,
\end{cases}
\end{equation}
with
\begin{align*}\label{notations}
u_k=(x_k,y_k,z_k)^\tau,\quad \mathcal{L}(u_k)=\Big(\mathcal{L}(x_k),\mathcal{L}(y_k),\mathcal{L}(z_k)\Big)^\tau,\quad u_{1,k}=\Big(u_k,\mathcal{L}(u_k),\mathcal{L}(u_{k-\tau}),\mathcal{L}(u_{k+\tau})\Big)^\tau.
\end{align*}
Here $\varphi_k\in C[-\tau,0]^n$, and define $f(0,\cdot)=0$, $y_k=z_k=0, k\in[-\tau,-1]$. So that
\begin{equation*}
\begin{aligned}
b,\sigma,f:\Omega\times\mathcal{T} \times \mathbb{R}^{3n} \times \left(\mathcal{P}_2(\mathbb{R}^n)\times\mathcal{P}_2(\mathbb{R}^n)\times\mathcal{P}_2(\mathbb{R}^{n})\right)^3\rightarrow \mathbb{R}^n.
\end{aligned}
\end{equation*}
Denote
\begin{align*}
\quad \mathcal{A}(k, u)
:=\Big((-f(k,u_{1,k}),b(k,u_{1,k}),\sigma(k,u_{1,k})\Big)^\tau.
\end{align*}
Then we state the following assumptions:
\begin{enumerate}
	\item[(H1)] There exists a constant $L>0$ independent of $\lambda,\tau$, such that for all processes $u=(x,y,z)$ and $u'=(x',y',z')$,
     {
    \begin{align*}
    \mathbb{E}\sum_{k=0}^{+\infty}e^{-\lambda k}\left(\|\mathcal{A}(k,u')-\mathcal{A}(k,u)\|^2\right)  
    \le L\mathbb{E}\sum_{k=0}^{+\infty}e^{-\lambda k}\|u'_k-u_k\|^2.
    \end{align*}
    }
For each $u$, $\mathcal{A}(\cdot,u)\in L_{\mathcal{F}}^{2,\lambda}(0,+\infty)$. Moreover, $\lambda$ is big enough as in Remark \ref{remark1}.

\item[(H2)] There exists a constant $L_1>0$, such that for all adapted processes $u'=(x',y',z')$, $u=(x,y,z)$ and $\widehat{u}=u'-u$,
\begin{align*}
      &\mathbb{E}\sum_{k=0}^{+\infty}e^{-\lambda k}\left[\|b(k,u'_{1,k})-b(k,u_{1,k})\|^2+\|\sigma(k,u'_{1,k})-\sigma(k,u_{1,k})\|^2\right]\\
        \le &L_1 \mathbb{E}\sum_{k=-\tau}^{+\infty}e^{-\lambda k}\|\widehat{x}_k\|^2+L_1 \mathbb{E}\sum_{k=0}^{+\infty}e^{-\lambda k}\|B_1^\tau\widehat{y}_k+B_2^\tau\mathbb{E}\widehat{y}_k+e^{-\lambda\tau}B_3^\tau\mathbb{E}\widehat{y}_{k+\tau}+D_1^\tau\widehat{z}_k+D_2^\tau\mathbb{E}\widehat{z}_k+e^{-\lambda\tau}D_3^\tau\mathbb{E}\widehat{z}_{k+\tau}\|^2,
\end{align*}
    for some given $\mathbb{R}^{n\times k}$-valued  {deterministic matrices} $B_i, D_i, i=1,2,3$.

\item[(H3)] There exists a constant $\beta>0$, such that
    \begin{align*}
&\mathbb{E}\sum_{k=0}^{+\infty}e^{-\lambda k}\langle\mathcal{A}(k,u')-\mathcal{A}(k,u),u'_k-u_k\rangle\\\le &-\beta\mathbb{E}\sum_{k=0}^{+\infty}e^{-\lambda k}\|B_1^\tau\widehat{y}_k+B_2^\tau\mathbb{E}\widehat{y}_k+e^{-\lambda\tau}B_3^\tau\mathbb{E}\widehat{y}_{k+\tau}+D_1^\tau\widehat{z}_k+D_2^\tau\mathbb{E}\widehat{z}_k+e^{-\lambda\tau}D_3^\tau\mathbb{E}\widehat{z}_{k+\tau}\|^2.
\end{align*}

\end{enumerate}

Then we show the existence and uniqueness of the solution to FBS$\Delta$Es (\ref{fbsde}).

\begin{theorem}\label{mainthe}
Assume that (H1)-(H3) hold. Then there exists a unique solution $u=(x,y,z)\in L_{\mathcal{F}}^{2,\lambda}\left(-\tau,+\infty ; \mathbb{R}^n\right)\times L_{\mathcal{F}}^{2,\lambda}\left(0,+\infty ; \mathbb{R}^n\right)\times L_{\mathcal{F}}^{2,\lambda}\left(0,+\infty ; \mathbb{R}^{n}\right)$ to the FBS$\Delta$E (\ref{fbsde}).
\end{theorem}

\begin{proof}
Now, we introduce the following FBS$\Delta$E with $\varepsilon\in[0,1]$, which becomes FBS$\Delta$E (\ref{fbsde}) when $\varepsilon=1$.
\begin{align}\label{fbsde1}
\begin{cases}\Delta x^\varepsilon_k=-(1-\varepsilon)\beta \left(B_1 u_{2,k}+B_2 \mathbb{E}u_{2,k}+B_3\mathbb{E} u_{2,k-\tau}\right)+\varepsilon b(k,u_{1,k})+\phi_k\\
\qquad\qquad+\left[-(1-\varepsilon)\beta \left(D_1 u_{2,k}+D_2 \mathbb{E}u_{2,k}+D_3\mathbb{E} u_{2,k-\tau}\right)+\varepsilon\sigma\left(k,u_{1,k}\right)+\psi_k\right]\Delta W_k, & k\in\mathcal{T}, \\ 
y_k+z_k\Delta W_k=e^{-\lambda}\left[y_{k+1}+\varepsilon f\left(k+1,u_{1,k+1}\right)+\gamma_{k+1}\right], & k\in\mathcal{T}, \\ x^\varepsilon_k=\varphi_k, & k \in[-\tau, 0], \\ 
\lim_{N\to+\infty}e^{-\lambda N}\mathbb{E}\left[\|y^\varepsilon_N\|^2+\|z^\varepsilon_N\|^2\right]=0.
\end{cases}
\end{align}
Here $\phi,\psi,\gamma\in L^{2,\lambda}_{\mathcal{F}}(0,+\infty;\mathbb{R}^n)$ are arbitrary given processes and 
\begin{align}\label{u2k}
u_{2,k}=\mathbf{1}_{\{k\ge 0\}}\left[B_1^\tau y_k+B_2^\tau\mathbb{E}y_k+e^{-\lambda\tau}B_3^\tau\mathbb{E}y_{k+\tau}+D_1^\tau z_k+D_2^\tau\mathbb{E}z_k+e^{-\lambda\tau}D_3^\tau\mathbb{E}z_{k+\tau}\right].
\end{align}

  {For the case $\varepsilon=0$, it is clear that the backward equation has a unique solution according to \cite{ji2024infinite}, and then the forward equation has a unique solution.
Our aim is to prove there exists  a unique solution to the case $\varepsilon=1$.} To show the  existence and uniqueness for the FBS$\Delta$E (\ref{fbsde1}) with $\varepsilon\in(0,1]$, we need the following Lemma.
\begin{lemma}\label{existence m>n}
Let (H1)-(H3) hold. If for some $\varepsilon_0\in[0,1)$, there exists a unique solution $(x^{\varepsilon_0},y^{\varepsilon_0},z^{\varepsilon_0})\in L_{\mathcal{F}}^{2,\lambda}\left(-\tau,+\infty ; \mathbb{R}^n\right)\times L_{\mathcal{F}}^{2,\lambda}\left(0,+\infty ; \mathbb{R}^n\right)\times L_{\mathcal{F}}^{2,\lambda}\left(0,+\infty ; \mathbb{R}^{n}\right)$ to the FBS$\Delta$Es (\ref{fbsde1})  {for any $\phi,\psi,\gamma$}. Then there exists a unique solution $(x^{\varepsilon_0+\delta},y^{\varepsilon_0+\delta},z^{\varepsilon_0+\delta})$ to FBS$\Delta$E (\ref{fbsde1}) for $\varepsilon=\varepsilon_0+\delta$  {for any $\phi,\psi,\gamma$}, where $\delta\in[0,\delta_0]$ for  some constant $\delta_0>0$ independent  {of} $\varepsilon_0$.
\end{lemma}
\begin{proof}
Let $u^i=(x^i,y^i,z^i), i=1,2$, be given processes in $L_{\mathcal{F}}^{2,\lambda}\left(-\tau,+\infty ; \mathbb{R}^n\right)\times L_{\mathcal{F}}^{2,\lambda}\left(0,+\infty ; \mathbb{R}^n\right)\times L_{\mathcal{F}}^{2,\lambda}\left(0,+\infty ; \mathbb{R}^{n}\right)$. Under the assumptions, we know there exists unique $U^i_t=(X^i_t,Y^i_t,Z^i_t),\, i=1,2$, which is the solution to the following FBS$\Delta$E,
\begin{align}\label{fbsde2}
\begin{cases}\Delta X^i_k=-(1-\varepsilon_0)\beta \left(B_1 U^i_{2,k}+B_2 \mathbb{E}U^i_{2,k}+B_3\mathbb{E} U^i_{2,k-\tau}\right)+\varepsilon_0 b(k,U^i_{1,k})+\phi_k
\\
\qquad\qquad+\delta\beta\left(B_1 u^i_{2,k}+B_2 \mathbb{E}u^i_{2,k}+B_3\mathbb{E} u^i_{2,k-\tau}\right)+\delta b(k,u^i_{1,k})
\\
\qquad\qquad+\left[-(1-\varepsilon_0)\beta \left(D_1 U^i_{2,k}+D_2 \mathbb{E}U^i_{2,k}+D_3\mathbb{E} U^i_{2,k-\tau}\right)+\varepsilon_0\sigma\left(k,U^i_{1,k}\right)+\psi_k\right]\Delta W_k,
\\ 
\qquad\qquad+\left[\delta\beta\left(D_1 u^i_{2,k}+D_2 \mathbb{E}u^i_{2,k}+D_3\mathbb{E} u^i_{2,k-\tau}\right)+\delta \sigma(k,u^i_{1,k})\right]\Delta W_k
\\
Y^i_k+Z^i_k\Delta W_k=e^{-\lambda}\Big[Y^i_{k+1}+\varepsilon_0 f\left(k+1,U^i_{1,k+1}\right)+\delta f\left(k+1,u^i_{1,k+1}\right)+\gamma_{k+1}\Big], & k\in\mathcal{T}, \\ 
X^i_k=\varphi_k, \qquad  {k \in[-\tau, 0]},
\\  
\lim_{N\to+\infty}e^{-\lambda N}\mathbb{E}\left[\|Y^i_N\|^2+\|Z^i_N\|^2\right]=0.
\end{cases}
\end{align}

Then we show the mapping $I: u\to U$ from $L_{\mathcal{F}}^{2,\lambda}\left(-\tau,+\infty ; \mathbb{R}^n\right)\times L_{\mathcal{F}}^{2,\lambda}\left(0,+\infty ; \mathbb{R}^n\right)\times L_{\mathcal{F}}^{2,\lambda}\left(0,+\infty ; \mathbb{R}^{n}\right)$ onto itself is a contraction under the norm
\begin{align*}
\|U\|^2_\lambda=\|(X,Y,Z)\|^2_\lambda=\mathbb{E}\left[\sum_{k=-\tau}^{+\infty}(1\land e^{-\lambda k})\|X_k\|^2+\sum_{k=0}^{+\infty}e^{-\lambda k}\left(\|Y_k\|^2+\|Z_k\|^2\right)\right].
\end{align*}
Denote $\widehat{U}_k=U^1_k-U^2_k, \,\widehat{u}_k=u^1_k-u^2_k$.  {Taking the difference of $\langle \widehat{X}_k, e^{-\lambda k}\widehat{Y}_k\rangle$ yields}
\begin{align*}
\Delta \langle \widehat{X}_k, e^{-\lambda k}\widehat{Y}_k\rangle=&\langle\Delta\widehat{X}_k, e^{-\lambda k}\widehat{Y}_k\rangle+\left\langle \widehat{X}_{k+1}, e^{-\lambda k}\left(e^{-\lambda}\widehat{Y}_{k+1}-\widehat{Y}_k\right)\right\rangle
\\
=&e^{-\lambda k}\Big\{-(1-\varepsilon_0)\beta\langle B_1 \widehat{U}_{2,k}+B_2 \mathbb{E}\widehat{U}_{2,k}+B_3\mathbb{E} \widehat{U}_{2,k-\tau},\widehat{Y}_k\rangle
\\
&\quad\qquad-(1-\varepsilon_0)\beta|\Delta W_k|^2\langle D_1 \widehat{U}_{2,k}+D_2 \mathbb{E}\widehat{U}_{2,k}+D_3\mathbb{E} \widehat{U}_{2,k-\tau},\widehat{Z}_k\rangle
\\
&\quad\qquad+\varepsilon_0\big[\langle b(k,U^1_{1,k})-b(k,U^2_{1,k}),\widehat{Y}_k\rangle+\langle \sigma(k,U^1_{1,k})-\sigma(k,U^2_{1,k}),\widehat{Z}_k\rangle
\\
&\qquad\qquad\qquad-e^{-\lambda}\langle f(k+1,U^1_{1,k+1})-f(k+1,U^2_{1,k+1}),\widehat{X}_{k+1}\rangle\big]\\
&\qquad\qquad+\delta\alpha(\widehat{U}_k,u^1_{1,k},u^2_{1,k})+M_k\Delta W_k\Big\},
\end{align*}
where $M_k$ is $\mathcal{F}_k$-measurable and $\alpha$ is the remainder term satisfies
\begin{align*}
\mathbb{E}\sum_{k=0}^{+\infty}e^{-\lambda k}\alpha(\widehat{U}_k,u^1_{1,k},u^2_{1,k})\le C\left(\|\widehat{U}\|_\lambda^2+\|\widehat{u}\|_\lambda^2\right).
\end{align*}
Here $C$ is a positive constant independent of $\lambda,\tau,\varepsilon_0,\delta$, which may vary from line to line. 
Through $\mathbb{E}\left[\Delta W_k|\mathcal{F}_k\right]=0$, $\mathbb{E}\left[\Delta W_k^2|\mathcal{F}_k\right]=1$, $\widehat{X}_0=0$ and $\lim_{N\to+\infty}\mathbb{E}\langle \widehat{X}_N, e^{-\lambda N}\widehat{Y}_N\rangle=0$, it follows
\begin{align}\label{444}
0=\mathbb{E}\sum_{k=0}^{+\infty}\Delta\langle \widehat{X}_k, e^{-\lambda k}\widehat{Y}_k\rangle\le \delta C\left(\|\widehat{U}\|^2_\lambda+\|\widehat{u}\|^2_\lambda\right)+\mathbb{E}\sum_{k=0}^{+\infty}e^{-\lambda k}\left[-(1-\varepsilon_0)\beta I_{1,k}+\varepsilon_0I_{2,k}\right],
\end{align}
where
\begin{align*}
&I_{1,k}=\langle B_1 \widehat{U}_{2,k}+B_2 \mathbb{E}\widehat{U}_{2,k}+B_3 \mathbb{E}\widehat{U}_{2,k-\tau},\widehat{Y}_k\rangle+\langle D_1 \widehat{U}_{2,k}+D_2 \mathbb{E}\widehat{U}_{2,k}+D_3\mathbb{E} \widehat{U}_{2,k-\tau},\widehat{Z}_k\rangle,
\\
&I_{2,k}=\langle b(k,U^1_{1,k})-b(k,U^2_{1,k}),\widehat{Y}_k\rangle+\langle \sigma(k,U^1_{1,k})-\sigma(k,U^2_{1,k}),\widehat{Z}_k\rangle
\\
&\qquad\qquad\quad-e^{-\lambda}\langle f(k+1,U^1_{1,k+1})-f(k+1,U^2_{1,k+1}),\widehat{X}_{k+1}\rangle.
\end{align*}
For $I_{1}$, we have
\begin{align}
\mathbb{E}\sum_{k=0}^{+\infty}e^{-\lambda k}I_{1,k}=&\mathbb{E}\sum_{k=0}^{+\infty}e^{-\lambda k}\left[\langle B_1 \widehat{U}_{2,k},\widehat{Y}_k\rangle+\langle B_2 \mathbb{E}\widehat{U}_{2,k},\widehat{Y}_k\rangle+\langle B_3 \mathbb{E}\widehat{U}_{2,k-\tau},\widehat{Y}_k\rangle\right]
\notag\\
&+\mathbb{E}\sum_{k=0}^{+\infty}e^{-\lambda k}\left[\langle D_1 \widehat{U}_{2,k},\widehat{Z}_k\rangle+\langle D_2 \mathbb{E}\widehat{U}_{2,k},\widehat{Z}_k\rangle+\langle D_3 \mathbb{E}\widehat{U}_{2,k-\tau},\widehat{Z}_k\rangle\right]
\notag\\
=&\mathbb{E}\sum_{k=0}^{+\infty}e^{-\lambda k}\left[\langle  \widehat{U}_{2,k},B_1^\tau\widehat{Y}_k\rangle+\langle \widehat{U}_{2,k},B_2^\tau\mathbb{E}\widehat{Y}_k\rangle+\langle \widehat{U}_{2,k},e^{-\lambda\tau}B_3^\tau\mathbb{E}\widehat{Y}_{k+\tau}\rangle\right]
\notag\\
&+\mathbb{E}\sum_{k=0}^{+\infty}e^{-\lambda k}\left[\langle  \widehat{U}_{2,k},D_1^\tau\widehat{Z}_k\rangle+\langle \widehat{U}_{2,k},D_2^\tau\mathbb{E}\widehat{Z}_k\rangle+\langle \widehat{U}_{2,k},e^{-\lambda\tau}D_3^\tau\mathbb{E}\widehat{Z}_{k+\tau}\rangle\right]
\notag\\
=&\mathbb{E}\sum_{k=0}^{+\infty}e^{-\lambda k}\|\widehat{U}_{2,k}\|^2.
\end{align}
For $I_2$, recall that $f(0,\cdot)=0$, we have
\begin{align}\label{666}
\mathbb{E}\sum_{k=0}^{+\infty}e^{-\lambda k}I_{2,k}=&\mathbb{E}\sum_{k=0}^{+\infty}e^{-\lambda k}\langle b(k,U^1_{1,k})-b(k,U^2_{1,k}),\widehat{Y}_k\rangle+\mathbb{E}\sum_{k=0}^{+\infty}e^{-\lambda k}\langle \sigma(k,U^1_{1,k})-\sigma(k,U^2_{1,k}),\widehat{Z}_k\rangle
\notag\\
&-\mathbb{E}\sum_{k=0}^{+\infty}e^{-\lambda k}\langle f(k,U^1_{1,k})-f(k,U^2_{1,k}),\widehat{X}_{k}\rangle
\notag\\
=&\mathbb{E}\sum_{k=0}^{+\infty}e^{-\lambda k}\langle\mathcal{A}(k,U^1)-\mathcal{A}(k,U^2),\widehat{U}_k\rangle
\notag\\
\le &-\beta\mathbb{E}\sum_{k=0}^{+\infty}e^{-\lambda k}\|\widehat{U}_{2,k}\|^2.
\end{align}
Combining with (\ref{444})-(\ref{666}), we conclude that
\begin{align*}
\mathbb{E}\sum_{k=0}^{+\infty}e^{-\lambda k}\|\widehat{U}_{2,k}\|^2\le\delta C\left(\|\widehat{U}\|^2_\lambda+\|\widehat{u}\|^2_\lambda\right).
\end{align*}
For the estimation of forward equation, through assumption (H2) and $\mathbb{E}\|\widehat{X}_k\|=0,\,k\in[-\tau,0]$, we have
\begin{align}\label{777}
(1-e^{-\lambda}C)\mathbb{E}\sum_{k=-\tau}^{+\infty}e^{-\lambda k}\|\widehat{X}_k\|^2\le \delta C\left(\|\widehat{U}\|^2_\lambda+\|\widehat{u}\|^2_\lambda\right)+C\mathbb{E}\sum_{k=0}^{+\infty}e^{-\lambda k}\|\widehat{U}_{2,k}\|^2
\le\delta C\left(\|\widehat{U}\|^2_\lambda+\|\widehat{u}\|^2_\lambda\right).
\end{align}
Then for backward equation, we have
\begin{align}\label{888}
(1-e^{-\lambda}C)\mathbb{E}\sum_{k=0}^{+\infty}e^{-\lambda k}\left(\|\widehat{Y}_k\|^2+\|\widehat{Z}_k\|^2\right)\le \delta C\left(\|\widehat{U}\|^2_\lambda+\|\widehat{u}\|^2_\lambda\right)+\mathbb{E}\sum_{k=0}^{+\infty}e^{-\lambda k}\|\widehat{X}_k\|^2\le \delta C\left(\|\widehat{U}\|^2_\lambda+\|\widehat{u}\|^2_\lambda\right).
\end{align}

\begin{remark}\label{remark1} 
$\lambda$ is big enough such that $1-e^{-\lambda}C>0$ in inequalities (\ref{777}) and (\ref{888}).
\end{remark}

Through inequalities (\ref{777}) and (\ref{888}), it follows
\begin{align*}
\|\widehat{U}\|^2_\lambda\le\delta C\left(\|\widehat{U}\|^2_\lambda+\|\widehat{u}\|^2_\lambda\right).
\end{align*}
Choose $\delta_0\le(3C)^{-1}$, it yields $\|\widehat{U}\|_\lambda^2\le\frac{1}{2}\|\widehat{u}\|^2_\lambda$. Thus, the mapping $I:u\to U$ is a contraction, so that there exists a unique fixed point in $L_{\mathcal{F}}^{2,\lambda}\left(-\tau,+\infty ; \mathbb{R}^n\right)\times L_{\mathcal{F}}^{2,\lambda}\left(0,+\infty ; \mathbb{R}^n\right)\times L_{\mathcal{F}}^{2,\lambda}\left(0,+\infty ; \mathbb{R}^{n}\right)$, which is the unique solution to FBS$\Delta$E (\ref{fbsde1}) for $\varepsilon=\varepsilon_0+\delta$.

This completes the proof of Lemma \ref{existence m>n}.
\end{proof}
Now, through induction, the existence and uniqueness of the solution to FBS$\Delta$E (\ref{fbsde1}) is proved for all $\varepsilon\in[0,1]$. So the proof of Theorem \ref{mainthe} is completed. 

\end{proof}

\section{Applications to optimal control}
In this section, we consider the following infinite horizon LQ McKean-Vlasov optimal control problem with law-delay:
\begin{align}\label{LQstate}
\left\{\begin{array}{ll}
\Delta x_k =A_1x_k+A_2\mathbb{E}x_k+A_3\mathbb{E}x_{k-\tau}+B_1v_k+B_2\mathbb{E}v_k+B_3\mathbb{E}v_{k-\tau},\\
\qquad\quad+[C_1x_k+C_2\mathbb{E}x_k+C_3\mathbb{E}x_{k-\tau}+D_1v_k+D_2\mathbb{E}v_k+D_3\mathbb{E}v_{k-\tau}]\Delta W_k&\quad k\in\mathcal{T},
\\x_k=\varphi_k, &k\in[-\tau,0],
\end{array}\right.
\end{align}
to minimize the cost functional
\begin{align}\label{LQcost}
    J(v)=\frac{1}{2}\mathbb{E}\sum_{k=0}^{+\infty}e^{-\lambda k}\left[\langle Q_1x_k,x_k\rangle+\langle Q_2\mathbb{E}x_k,\mathbb{E}x_k\rangle+\langle Rv_k,v_k\rangle\right].
\end{align}

Here $Q_1, Q_2\ge0$ are $\mathcal F_0$-measurable bounded random symmetric matrices and $R>0$ is a deterministic symmetric matrix such that $\lambda_{min}(R)\ge \alpha_1>0$, and $A_i,B_i,C_i,D_i,i=1,2,3,$ are deterministic matrices with suitable dimensions. 

Denote by $\mathbb{U}$ the set of admissible control process
$v=(v_k)_{k\ge-\tau}$ taking values in $\mathbb{R}^k$ such that $v_k=\theta_k, k\in[-\tau,-1]$ and $\mathbb{E}\sum_{k=-\tau}^{+\infty}(1\land e^{-\lambda k})\|v_k\|^2 <+\infty$. Here  {$\theta_k$} is a given $\mathcal{F}_0$-measurable process.
Our aim is to find the optimal control $v^*$ such that
\begin{align*}
  {   J(v^*)=\inf_{v\in \mathbb{U}}J(v)}
.\end{align*}

\begin{theorem}
The function given by
\begin{align}\label{v}
v_k=\begin{cases}\theta_k, & k \in[-\tau, -1], 
\\
 {-R^{-1}\left(B_1^\tau y_k+B_2^\tau \mathbb{E}y_k+e^{-\lambda \tau}B_3^\tau \mathbb{E}y_{k+\tau}+D_1^\tau z_k+D^\tau _2\mathbb{E}z_k+e^{-\lambda \tau}D^\tau _3\mathbb{E}z_{k+\tau}\right),} & k \in\mathcal{T},
\end{cases}
\end{align}    
is the unique optimal control for LQ problem (\ref{LQstate}) and (\ref{LQcost}), where $(x,y,z)$ is the solution to the following FBS$\Delta$E:
\begin{equation}\label{fbsde7}
\begin{cases}
\Delta x_k =A_1x_k+A_2\mathbb{E}x_k+A_3\mathbb{E}x_{k-\tau}+B_1v_k+B_2\mathbb{E}v_k+B_3\mathbb{E}v_{k-\tau}\\
\qquad\quad+[C_1x_k+C_2\mathbb{E}x_k+C_3\mathbb{E}x_{k-\tau}+D_1v_k+D_2\mathbb{E}v_k+D_3\mathbb{E}v_{k-\tau}]\Delta W_k,&\quad k\in\mathcal{T},
\\
y_k= {e^{-\lambda}\big[y_{k+1}+A_1^\tau y_{k+1}+ A_2^\tau\mathbb{E}[ y_{k+1}]+ e^{-\lambda \tau}A_3^\tau\mathbb{E}[ y_{k+\tau+1}]+Q_1x_{k+1}+\mathbb{E}Q_2\mathbb{E}x_{k+1}}\\
\qquad\qquad\qquad {+C_1^\tau z_{k+1}+ C_2^\tau \mathbb{E}[z_{k+1}]+ e^{-\lambda \tau}C_3^\tau \mathbb{E}[z_{k+\tau+1}]\big]}-z_k\Delta W_k, & k \in\mathcal{T}, \\
x_k=\varphi_k, &k\in[-\tau,0],\\ 
 \lim_{N\to+\infty}e^{-\lambda N}\mathbb{E}\left[\|y_N\|^2+\|z_N\|^2\right]=0,\end{cases}
\end{equation}
where $v_k$ is given by  equation (\ref{v}).
\end{theorem}

\begin{proof}
We first show the existence and uniqueness of the solution to FBS$\Delta$E (\ref{fbsde7}). It is clear that FBS$\Delta$E (\ref{fbsde7}) satisfies (H1) and (H2). 
For (H3), it is not difficult to obtain
{
\begin{align*}
&\mathbb{E}\sum_{k=0}^{+\infty}e^{-\lambda k}\langle\mathcal{A}(k,u')-\mathcal{A}(k,u),u'_k-u_k\rangle\\
=&-\mathbb{E}\sum_{k=0}^{+\infty}e^{-\lambda k}\left[\langle Q_1\widehat{x}_k,\widehat{x}_k\rangle+\langle Q_2\mathbb{E}\widehat{x}_k,\mathbb{E}\widehat{x}_k\rangle\right]+\mathbb{E}\sum_{k=0}^{+\infty}e^{-\lambda k}\langle \widehat{v}_k,\widehat{u}_{2,k}\rangle\\
\le&-\mathbb{E}\sum_{k=0}^{+\infty}e^{-\lambda k}\langle R^{-1}\widehat{u}_{2,k},\widehat{u}_{2,k}\rangle.
\end{align*}
}
Choosing $\beta=\left(
\lambda_{\max}(R)
\right)^{-1}$, we show (H3) holds.

Now we show the $v$ given by  equation (\ref{v}) is an optimal control. Let  {$v'$} be any other admissible control and  {$x'$} be the corresponding state process. Notice that $\langle \widehat{x}_0,y_0\rangle=\lim_{N\to+\infty}e^{-\lambda N}\langle\widehat{x}_N,y_N\rangle=0 $, so
\begin{align*}
0=&\mathbb{E}\sum_{k=0}^{+\infty}\Delta \langle \widehat{x}_k,e^{-\lambda k}y_k\rangle\\
    =&-\mathbb{E}\sum_{k=0}^{+\infty}e^{-\lambda k}\left[\langle Q_1x_k,\widehat{x}_k\rangle+\langle Q_2\mathbb{E}x_k,\mathbb{E}\widehat{x}_k\rangle+\langle Rv_k,\widehat{v}_k\rangle\right]
    \\
    \ge&-\frac{1}{2}\mathbb{E}\sum_{k=0}^{+\infty}e^{-\lambda k}\left[\langle Q_1x'_k,x_k'\rangle-\langle Q_1x_k,x_k\rangle+\langle Q_2\mathbb{E}x'_k,\mathbb{E}x_k'\rangle-\langle Q_2\mathbb{E}x_k,\mathbb{E}x_k\rangle+\langle Rv'_k,v'_k\rangle-\langle Rv_k,v_k\rangle\right].
\end{align*}
The above inequality means
\begin{align*}
J(v')-J(v)\ge0.
\end{align*}
For uniqueness,  assume that  both $v'$ and $v$ are optimal controls, $x'$ and $x$ are corresponding state processes, respectively. It is easy to see that $\frac{x'+x}{2}$ is the corresponding state process to $\frac{v'+v}{2}$.{ Recall that $\lambda_{min}(R)\ge\alpha_1>0,\, a.s.,$ and assume that there exists $ \alpha_2\ge 0$, such that $ J(v')=J(v)=\alpha_2.$}
So we have
\begin{align*}
2\alpha_2=&J(v')+J(v)\\
           =&\frac{1}{2}\mathbb{E}\sum_{k=0}^{+\infty}e^{-\lambda k}[\langle Q_1x_k',x_k'\rangle+\langle Q_1 x_k,x_k\rangle+\langle Q_2\mathbb{E}x_k',\mathbb{E}x_k'\rangle+\langle Q_2\mathbb{E}x_k,\mathbb{E}x_k\rangle+\langle Rv_k',v_k'\rangle+\langle Rv_k,v_k\rangle]
           \\
           \ge&\mathbb{E}\sum_{k=0}^{+\infty}e^{-\lambda k} \left[\left\langle Q_1\frac{x_k'+x_k}{2},\frac{x_k'+x_k}{2}\right\rangle+\left\langle Q_2\frac{\mathbb{E}x_k'+\mathbb{E}x_k}{2},\frac{\mathbb{E}x_k'+\mathbb{E}x_k}{2}\right\rangle+\left\langle R\frac{v_k'+v_k}{2},\frac{v_k'+v_k}{2}\right\rangle\right]
           \\
           &+\mathbb{E}\sum_{k=0}^{+\infty}e^{-\lambda k}\left\langle R\frac{v_k'-v_k}{2},\frac{v_k'-v_k}{2}\right\rangle\\
           =&2J\Big(\frac{v'+v}{2}\Big)+\mathbb{E}\sum_{k=0}^{+\infty}e^{-\lambda k}\left\langle R\frac{v_k'-v_k}{2},\frac{v_k'-v_k}{2}\right\rangle \\
           \ge&2\alpha_2+\frac{\alpha_1}{4}\mathbb{E}\sum_{k=0}^{+\infty}e^{-\lambda k}|v_k'-v_k|^2,
\end{align*}
which implies that $\mathbb{E}\sum_{k=0}^{+\infty}e^{-\lambda k}|v_k'-v_k|^2=0$. The proof completes.
\end{proof}


\bibliographystyle{cas-model2-names}

\bibliography{cas-refs}

@article{andersson2011maximum,
  title={A maximum principle for {SDEs} of mean-field type},
  author={Andersson, Daniel and Djehiche, Boualem},
  journal={Applied Mathematics \& Optimization},
  volume={63},
  number={3},
  pages={341--356},
  year={2011},
  publisher={Springer}
}

@article{buckdahn2016stochastic,
  title={A stochastic maximum principle for general mean-field systems},
  author={Buckdahn, Rainer and Li, Juan and Ma, Jin},
  journal={Applied Mathematics \& Optimization},
  volume={74},
  number={3},
  pages={507--534},
  year={2016},
  publisher={Springer}
}

@article{nie2022extended,
  title={Extended mean-field control problem with partial observation},
  author={Nie, Tianyang and Yan, Ke},
  journal={ESAIM: Control, Optimisation and Calculus of Variations},
  volume={28},
  pages={17},
  year={2022},
  publisher={EDP Sciences}
}

@article{dong2022maximum,
  title={Maximum principle for discrete-time stochastic control problem of mean-field type},
  author={Dong, Bozhang and Nie, Tianyang and Wu, Zhen},
  journal={Automatica},
  volume={144},
  pages={110497},
  year={2022},
  publisher={Elsevier}
}

@article{buckdahn2026global,
  title={A global stochastic maximum principle for mean-field forward-backward stochastic control systems with quadratic generators},
  author={Buckdahn, Rainer and Li, Juan and Li, Yanwei and Wang, Yi},
  journal={The Annals of Applied Probability},
  volume={36},
  number={1},
  pages={275--318},
  year={2026},
  publisher={Institute of Mathematical Statistics}
}

@article{chen2019maximum,
  title={Maximum principle for delayed stochastic mean-field control problem with state constraint},
  author={Chen, Li and Wang, Jiandong},
  journal={Advances in Difference Equations},
  volume={2019},
  number={1},
  pages={348},
  year={2019},
  publisher={Springer}
}

@article{zhang2021stochastic,
  title={Stochastic maximum principle of mean-field jump--diffusion systems with mixed delays},
  author={Zhang, Feng},
  journal={Systems \& Control Letters},
  volume={149},
  pages={104874},
  year={2021},
  publisher={Elsevier}
}

@article{guo2024stochastic,
  title={Stochastic maximum principle for generalized mean-field delay control problem},
  author={Guo, Hancheng and Xiong, Jie and Zheng, Jiayu},
  journal={Journal of Optimization Theory and Applications},
  volume={201},
  number={1},
  pages={352--377},
  year={2024},
  publisher={Springer}
}

@article{lin2014maximum,
  title={A maximum principle for optimal control of discrete-time stochastic systems with multiplicative noise},
  author={Lin, Xiangyun and Zhang, Weihai},
  journal={IEEE Transactions on Automatic Control},
  volume={60},
  number={4},
  pages={1121--1126},
  year={2014},
  publisher={IEEE}
}

@article{dong2023maximum,
  title={A maximum principle for discrete-time stochastic optimal control problem with delay},
  author={Dong, Bozhang and Nie, Tianyang and Wu, Zhen},
  journal={Systems \& Control Letters},
  volume={181},
  pages={105631},
  year={2023},
  publisher={Elsevier}
}

@article{ji2024infinite,
  title={Infinite Horizon Backward Stochastic Difference Equations and Related Stochastic Recursive Control Problems},
  author={Ji, Yuanyuan and Zhang, Qi},
  journal={SIAM Journal on Control and Optimization},
  volume={62},
  number={5},
  pages={2750--2775},
  year={2024},
  publisher={SIAM}
}

@article{wu2022maximum,
  title={Maximum principle for discrete-time stochastic optimal control problem and stochastic game.},
  author={Wu, Zhen and Zhang, Feng},
  journal={Mathematical Control \& Related Fields},
  volume={12},
  number={2},
  year={2022}
}

@article{niu2026stochastic,
  title={Stochastic maximum principle for fully coupled nonlinear {FBS$\Delta$Es} under generalized monotonicity and {LQ} control applications},
  author={Niu, Zhipeng and Moon, Jun and Meng, Qingxin},
  journal={Systems \& Control Letters},
  volume={208},
  pages={106328},
  year={2026},
  publisher={Elsevier}
}

@article{hu1995solution,
  title={Solution of forward-backward stochastic differential equations},
  author={Hu, Ying and Peng, Shige},
  journal={Probability Theory and Related Fields},
  volume={103},
  pages={273--283},
  year={1995},
  publisher={Springer}
}

@article{peng1999fully,
  title={Fully coupled forward-backward stochastic differential equations and applications to optimal control},
  author={Peng, Shige and Wu, Zhen},
  journal={SIAM Journal on Control and Optimization},
  volume={37},
  number={3},
  pages={825--843},
  year={1999},
  publisher={SIAM}
}

@article{ma1994solving,
  title={Solving forward-backward stochastic differential equations explicitly—a four step scheme},
  author={Ma, Jin and Protter, Philip and Yong, Jiongmin},
  journal={Probability Theory and Related Fields},
  volume={98},
  number={3},
  pages={339--359},
  year={1994},
  publisher={Springer}
}

@article{song2023solvability,
  title={Solvability of general fully coupled forward--backward stochastic difference equations with delay and applications},
  author={Song, Teng},
  journal={Optimal Control Applications and Methods},
  volume={44},
  number={6},
  pages={3011--3032},
  year={2023},
  publisher={Wiley Online Library}
}



\end{document}